\documentclass[12pt]{article}
\usepackage{graphicx,psfrag}

\begin{document}
 
\begin{center}
               {\bf{How to beat Capablanca}} \\ 
               {\em{by Thotsaporn Thanatipanonda}} \\
               {\em{Email: thot@math.rutgers.edu}} \\
\end{center}

{\bf{   1. Introduction}}\\

      The first thing Capablanca mentions in  his book, 
$\underline{Chess}$  $\underline{Fundamentals}$, is how to checkmate with 
rook, as in the diagram below.\\

\begin{figure}[h]
  \begin{center}
\resizebox{7cm}{7cm}{\includegraphics{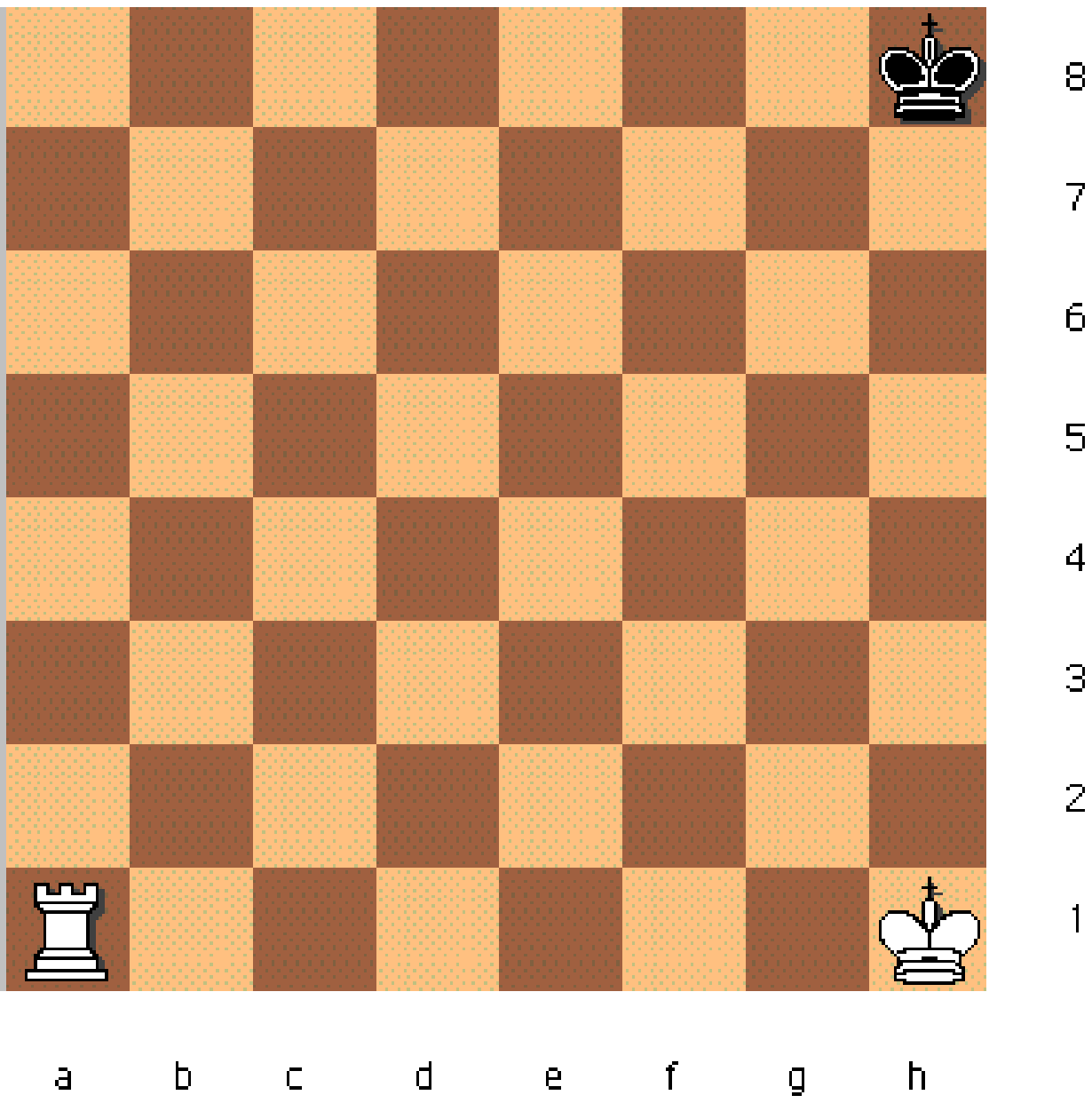}}

  \end{center}
\end{figure}

\noindent       Capablanca writes, "In this position the power of the rook 
is 
demonstrated by the first move, Ra7, which immediately confines the black king to the last rank, ...". 
Capablanca did not give the fastest way to checkmate. 
With {\bf{1.Ra7}}, the fastest way for White to checkmate Black king is 10 
moves. 
While with {\bf{1.Rg1}}, White can force a checkmate in 9 moves.\\

{\bf{2. On an $m \times n$ board}} \\
    
\noindent    For the general $m \times n$ board with White king (WK) on 
$(m,1)$ , White rook (WR) on (1,1), and Black king (BK) on $(m,n)$ where 
$m\geq4$ and $n\geq5$, we 
define $FM(m,n)$ to be the smallest number of moves for white to checkmate 
Black king.  \\

\noindent {\bf{Theorem}}:$FM(m,n)$ =  $\left \{ \begin{array}{ll} 
                 n & \mbox{if n is odd};\\
               n+1 & \mbox{if n is even}. 
                  \end{array}
                 \right.$ \\

\noindent Let's call the right side $G(m,n)$. We will prove this by \\
1) Showing a sequence of White moves that force the checkmate in 
$G(m,n)$ moves.\\
2) Showing that Black has a strategy to survive up to $G(m,n)-1$ moves.\\

\noindent {\bf{Lemma 1}} In the given postion on an $m \times n$ board, 
White 
can give a checkmate in $n$ moves if $n$ is odd. \\

I will give the sequence of white moves so that for all the choices of 
Black moves, white can give the checkmate in $n$ moves.\\

First move: 1.WR(m-1,1)  \hspace*{2.0cm}      BK(m,n-1) (only move)\\

Then White can make a sequence of White king moves up one square 
no matter what Black responses are. Black's only responses are moving the 
king 
up and 
down along the $m^{th}$column. But this can not interrupt White king because of the parity. 
It will take $n-4$ moves for white king to move up to the square 
$(m,n-3)$. We have the diagram below with White to move.

\begin{figure}[h]
  \begin{center}
\resizebox{7cm}{7cm}{\includegraphics{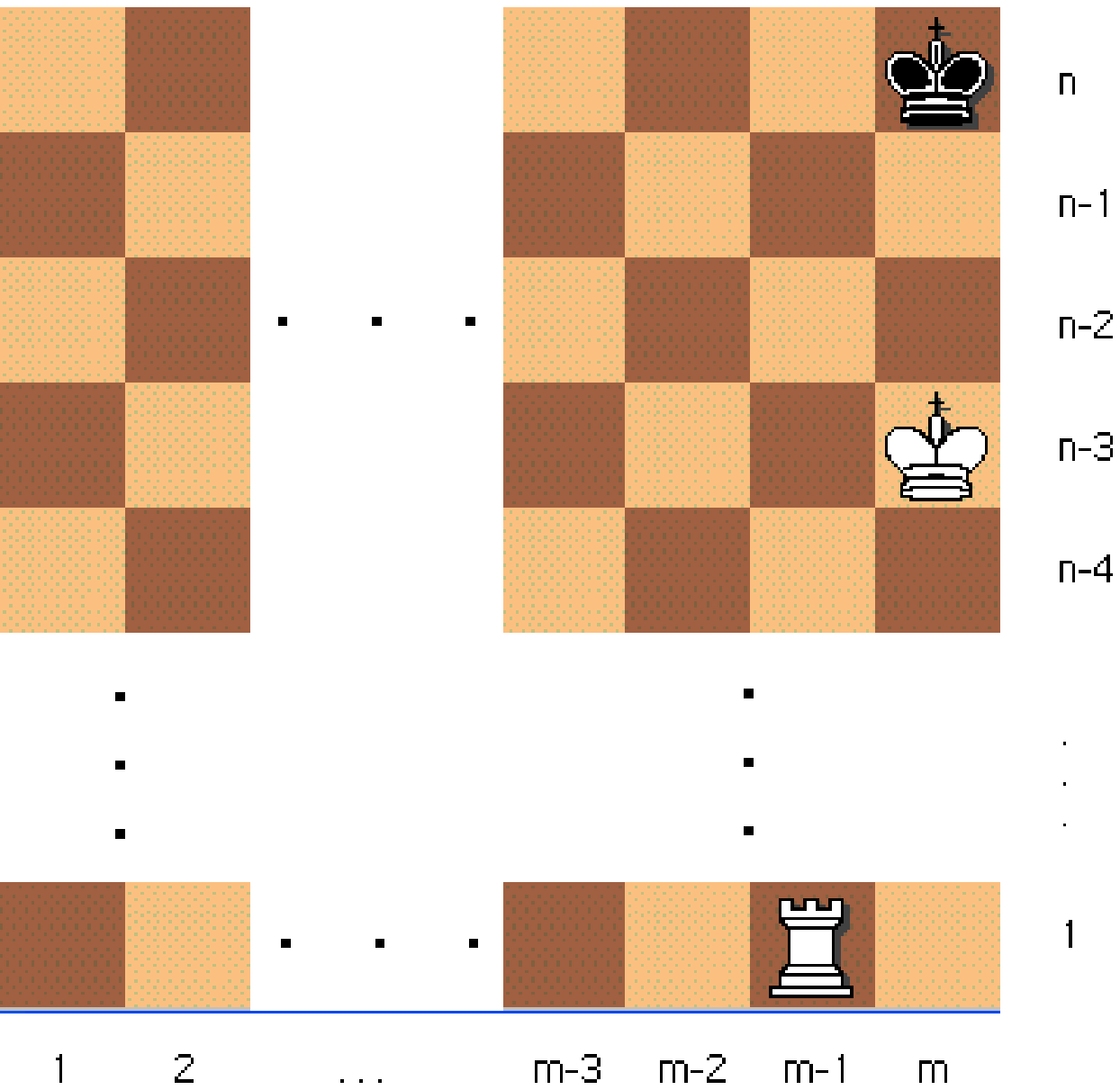}}

  \end{center}
\end{figure}
      
\noindent White can finish this off by playing the forcing move\\

\noindent 1. WK(m-1,n-2)  \hspace*{1.57cm} BK(m-1,n) (only move)\\
2. WR(m-2,1)   \hspace*{2.0cm}  BK(m,n)   (only move)\\
3. WR(m-2,n)   \hspace*{2.0cm}  Checkmate \\

\noindent The total number of moves is $1 + (n-4) + 3 = n$. \\

\noindent {\bf{Lemma 2}} White can force a checkmate in $n+1$ moves, if 
$n$ is even. \\

The strategy is almost the same as in Lemma 1. However, because of the 
parity we have to make one waiting move, i.e.,WR(m-1,2) while White king 
tries to move upward. Since $n\geq 5$, we have enough room for this plan.\\ 

\noindent {\bf{Note}}: The first move WR(m-1,1)
with the  strategy above gives White the fastest way to mate. As we can see it requires 
at least $n-3$ (resp. $n-2$ moves) for $n$ odd (resp. even) for White King to 
move up the board and at least 2 more Rook moves before White can force a checkmate.\\

\noindent {\bf{Lemma 3}} Black has a strategy to survive up to $n-1$ moves, if $n$ is 
odd. \\

\noindent The general plan for Black to survive is to move to the 
middle of the board as much as possible. At some point White has to move 
his Rook to restrict the possible moves of Black king. Then he must use 
White king to push Black king to the edge of the board. In addition, the 
sooner the Rook moves, the better. Therefore we will assume the first move is a 
Rook move. Note that the only mating positions are when Black king is 
at the edge of the board. Now we consider 2 cases of the first Rook 
move:\\

\noindent {\bf{Case 1}} The first rook move is vertical (moving along 
a column).\\  
        
In this case, Black tries to move to the middle of the board as much as 
possible. Once Black king gets there, he will try to stay there as long 
as possible before he is forced to the corner. Below is an example where 
the size of the board is 
$9\times9$.\\

\noindent 1. WR(1,8)  \hspace*{2.0cm}     BK(8,9)\\  
2. WK(8,2)  \hspace*{2.0cm}     BK(7,9)\\
3. WK(7,3)  \hspace*{2.0cm}     BK(6,9)\\
4. WK(6,4)  \hspace*{2.0cm}     BK(5,9)\\
5. WK(5,5)  \hspace*{2.0cm}     BK(6,9)\\
6. WK(5,6)  \hspace*{2.0cm}     BK(5,9)\\
7. WK(4,7)  \hspace*{2.0cm}     BK(6,9)\\
8. WK(5,7)  \hspace*{2.0cm}     BK(7,9)\\
9. WK(6,7)  \hspace*{2.0cm}     BK(8,9)\\
10.WK(7,7)  \hspace*{2.0cm}     BK(9,9)\\
11.WK(8,7)  \hspace*{2.0cm}     BK(8,9)\\
12.WR(1,9)  \hspace*{2.0cm}     Checkmate  \\

On $m \times n$ board, we can see that White has to move his Rook twice, and
White king 
moves up the board in $n-3$ moves and chases Black king back to the corner with at 
least another $\lceil(\frac{m}{2})\rceil-1$ moves. For example on a $9\times9$ 
board, it takes $2+(9-3)+(5-1) = 12$ moves. This is the best that White can do.
Therefore $F(m,n) = 2+(n-3)+(\lceil(\frac{m}{2})\rceil-1) \geq n $ moves 
(since  $m \geq 4$ 
).\\

\noindent {\bf{Case 2}} The first Rook move is horizontal (moving along a 
row). \\
(We exclude the move $R(m-1,1)$ since we already know that in that case the 
fastest number of moves to mate is n.)\\

Black king will try to move down and toward the middle of the 
board as much as 
possible. Once he gets blocked by the rook and the White king, Black king moves along the 
row (otherwise Black King will just try to stay in the middle of the 
board). He could move along the row since the Rook is not at column $m-1$. 
The example below illustrates the situation:\\

\newpage
\begin{figure}[h]
  \begin{center}
\resizebox{8cm}{8cm}{\includegraphics{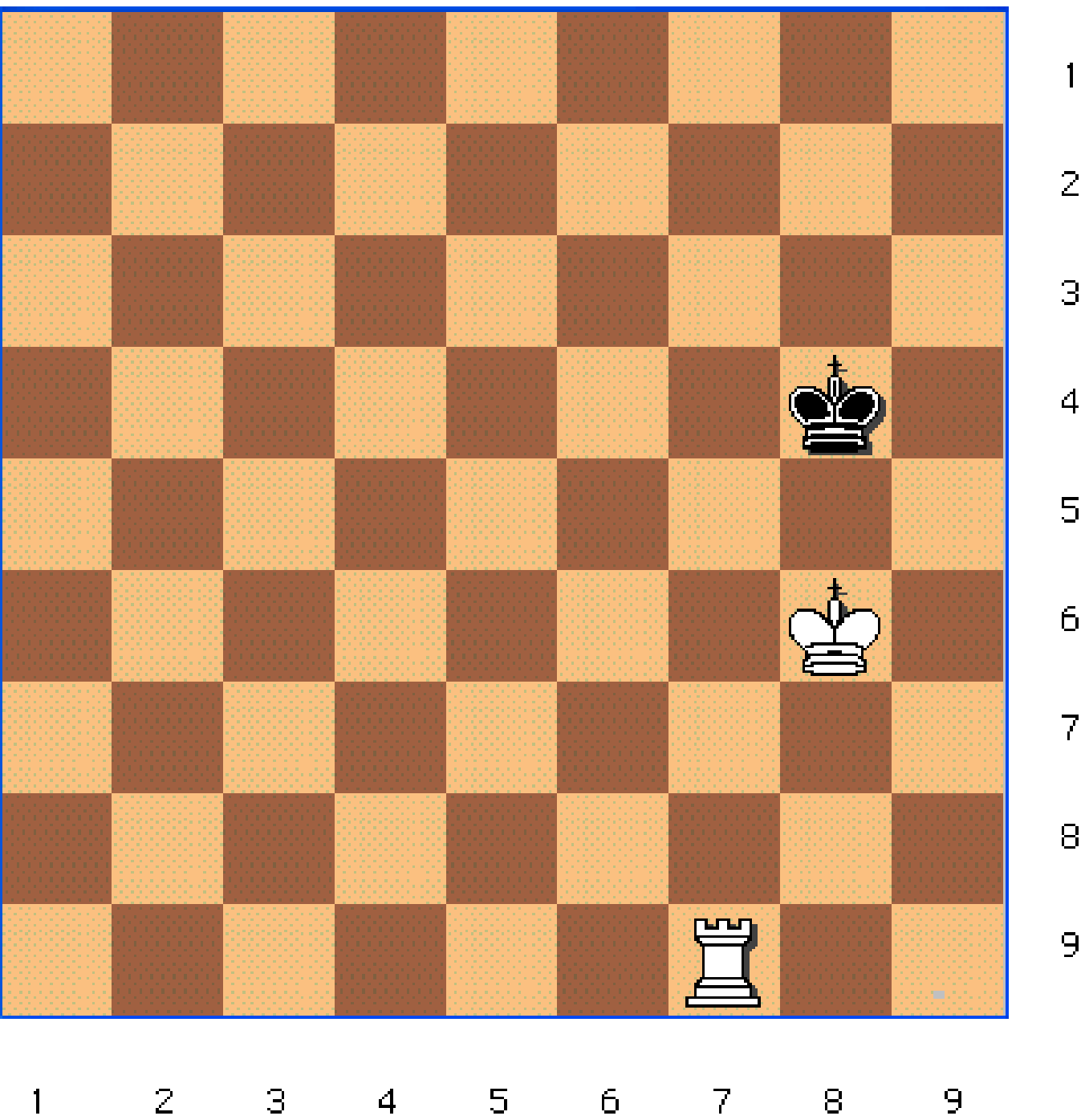}}
                 
  \end{center}
\end{figure}

1. ...   BK(9,6)   Now white has two choices: 2.WK(7,5) or 2.WR(7,6).
For 2.WK(7,5), Black's response will be BK(8,7) and the best white response 
is 3.WK(8,5). In this case, White King needs one extra move to get to the row 
$n-3$.\\

Overall, White has to move the rook twice (the first rook move and 
the checkmate move); the king takes $n-3$ moves to go up the board and at least an 
extra  
move from king or rook of the two choices above. This makes the least possible number of moves to 
checkmate at least $ \geq 2+(n-3)+1 = n.$\\  
 
\noindent {\bf{Lemma 4}} Black has a strategy to survive up to $n$ moves, when $n$ is 
even.\\ 

The general plan in this case is similar. The strategy in Case 1 also works here. In 
Case 2, Black now has an advantage of parity. This will give him an extra move to 
survive 
when he faces White king in the middle of the board. The details of the proof are left as an 
exercise to the reader.\\

By Lemmas 1, 2, 3 and 4, we conclude the Theorem. \\

\noindent {\bf{3. About Jos\'{e} Raul Capablanca}} \\

\noindent  Jose Raul Capablanca was born in Havana, Cuba on November 19, 
1888. He is regarded as one of the most gifted chess players of all time. According to 
Capablanca, at the age of four, he learned the rules of chess by watching 
his father 
play chess with a friend. In 1921, he won the world champion title from  
mathematician and chess player Emanuel Lasker, who held the title for 27 
years. Capablanca 
was the third official world chess champion. He lost his title to Alexander 
Alekhine in 1927 and died in 1942.\\

\noindent  The move that Capablanca suggested is Case 1 of Lemma 4.
Black could survive at least $2+(n-3)+(\lceil(\frac{m}{2})\rceil-1) = 
2+5+3 = 1$\\

\noindent  I want to thank my advisor, Doron Zeilberger for suggesting 
the problem and constantly giving help and support.\\

\noindent {\bf{Reference}}: {\em{Chess Fundamentals}}, Jos\'{e} Raul 
Capablanca; Everyman Chess, 1999. First published in 1921. \\ 

\end{document}